\documentclass[11pt]{amsart}

\usepackage[utf8]{inputenc}
\usepackage[T1]{fontenc}
\usepackage{amsmath,amssymb,amsthm,mathtools}
\usepackage{enumitem}
\usepackage{hyperref}
\usepackage[margin=1in]{geometry}

\hypersetup{colorlinks=true,linkcolor=blue,citecolor=blue,urlcolor=blue}

\newtheorem{theorem}{Theorem}[section]

\newtheorem{proposition}[theorem]{Proposition}
\newtheorem{corollary}[theorem]{Corollary}

\theoremstyle{definition}

\newtheorem{remark}[theorem]{Remark}

\DeclareMathOperator{\Exp}{Exp}
\DeclareMathOperator{\Var}{MSD}

\DeclareMathOperator{\Ree}{Re}
\newcommand{\E}{\mathbb{E}}
\newcommand{\Pp}{\mathbb{P}}
\newcommand{\R}{\mathbb{R}}

\newcommand{\dd}{\mathrm{d}}
\newcommand{\eqd}{\stackrel{d}{=}}

\newcommand{\ol}{\overline}
\newcommand{\ul}{\underline}
\newcommand{\RV}{\mathrm{RV}}

\title{Steady-state of L\'evy process with general resetting: unreflected and reflected cases}
\author{Zbigniew Palmowski}
\address{Faculty of Pure and Applied Mathematics,
	Wroc\l aw University of Science and Technology
	Wyb. Wyspia\'nskiego 27, 50-370 Wroc\l aw, Poland}
\email{zbigniew.palmowski@pwr.edu.pl}
\thanks{Zbigniew Palmowski was partially supported by the National Science Centre (Poland): grant 2023/51/B/ST1/01270.}

\author{Maria Vlasiou}
\address{University of Twente, the Netherlands}
\email{m.vlasiou@utwente.nl}
\date{\today}

\begin{document}
\maketitle
\begin{abstract}
We study (unreflected and reflected) L\'evy processes that are modified at independent Poisson review epochs by a random state-dependent reset map.  The first model is unreflected: between resettings the level evolves as a L\'evy process, and at a resetting 
epoch the new level is $F$ of the level seen just before the resetting.
The second model is the reflected L\'evy process at $0$, but with reset maps chosen from a wider transform-solvable class.  In this work, we give contraction-type ergodicity criteria for the existence of the steady-state.
We identify
this steady-state of the system
as the solution of certain stochastic fixed-point equation.
This allows to find the Laplace-transform, moments recursions, and regularly-varying-tail formulae.
For non-reflected case, the examples  include affine deterministic resets $F(x)=cx+b$, random-coefficient resets $F(x)=Cx+b$ with two-point $C$ for the finite-activity L\'evy process, spectrally positive 
L\'evy flight, reflected Brownian process.  For the reflected cased we analyse affine resetting mechanism $F(x)=cx$, the alternating resetting $F(x)=(B-x)^+$, and uniform proportional resetting $F(x)=Cx$ with $C$ having uniform distribution.
\end{abstract}

\section{Introduction}

There has been a vast amount of papers related to the resetting mechanism.
This mechanisms arise in search processes (see, e.g., \cite{Bel, Evans, White}), physics (see, e.g., \cite{18}), biology (see, e.g., \cite{4}), queueing theory (see, e.g., \cite{MR2023017}) and in many other
research fields; see also the survey~\cite{Gupta}.
The usual resetting mechanism is based either by sending the process to $0$ (or any other fixed state) or by multiplying by possibly random amplitude at independent Poisson epochs.
Lately, there have been efforts to consider more general resetting mechanisms, see for example \cite{Martinetal} and references therein.
Based on 
\cite{PalmowskiVlasiou2011}, we also consider more general possibly random resetting mechanisms $F(x)$, where $x$ is state prior resetting,
that is applied to L\'evy processes $X$ 
and its reflected version $Y$ at $0$.

For unreflected L\'evy process $X$, we give a contraction ergodicity criterion for the resetting mechanism to get steady-state of our system.
For the stationary law we give a general transform equation and moment recursions for arbitrary reset maps $F$ whenever powers of $F(Z)$ can be expressed through the moments of $Z$.
For the resetting mechanism $F(x)=cx+b$ (where $c,b$ are constants), using a perpetuity representation, we obtain an infinite-product transform of the stationary law. We also identify 
the moment recursions.
Finally, we also present the behaviour of the asymptotic tail of the law of this stationary distribution under regular variation of the positive jumps.
We managed also to analyse another interesting case when $F(x)=Cx+b$ with two-point and three-point distribution of $C$.

For the reflection at $0$ case we consider L\'evy process without negative jumps $X$. Such process is called spectrally positive.
For spectrally positive L\'evy process $Y$ reflected at $0$ we also give sufficient conditions for the existence of the steady-state,
and wrire an appropriate recursion identifying the limiting law.
We analyse in detail the retting maps $F(x)=cx$, $F(x)=(B-x)^+$ with exponential $B$, and $F(x)=Cx$ with $C$ having uniform distribution $U(0,1)$.
We also include the tail of distrubution of the steady-state in the case when jumps have
regularly varying jumps.

The present paper has two aims. Firstly, for very general process between resettings and quite general resetting mechanisms, we give the ergodicity conditions and
show that the limiting stationary distributions satisfies a one-step recursive equation. This equations codes some recursive procedures to find moments of
this stationary distribution and allows to identify the tail behaviour when jumps of the process between resettings is regularly varying.
The second aim is to analyse the process reflected at $0$ and then reset it. We show that, to some extend, the same ideas as for the process without reflection
will still work.

The paper is organised as follows. In the next section we analyse the case without reflection and then we focus on the reflected case.
%
%
Throughout this paper random reset coefficients are denoted by capital letters, such as $C$; deterministic reset constants are denoted by lower-case letters, such as $c$ and $b$.

\section{The unreflected reset model}\label{sec:unreflected}
\subsection{Notations}
For a real-valued L\'evy process $X$, write
\[
        \psi(s):=\log < e^{-sX(1)}>,
\]
whenever the expectation is finite.  Then, for an independent exponential random variable  $e_q\sim\Exp(q)$ with intensity $q$, we denote
\begin{equation}\label{eq:unref-exp-transform}
        m(s):=< e^{-sX(e_q)}>=\int_0^\infty q e^{-qt}e^{t\psi(s)}\,\dd t
        =\frac{q}{q-\psi(s)},
\end{equation}
for $s$ in the domain
\[
        \mathcal D_q:=\{s:\Ree(q-\psi(s))>0\}.
\]
In the reflected spectrally positive part of the paper, we use the standard fluctuation notation
\[
        \ol X(t):=\sup_{0\leq u\leq t}X(u),\qquad
        \ul X(t):=\inf_{0\leq u\leq t}X(u).
\]
Let $\Phi(q)$ be the largest non-negative solution of
\begin{equation}\label{eq:Phi-root}
        \psi(\Phi(q))=q.
\end{equation}
For a spectrally positive L\'evy process, 
the Wiener--Hopf factorisation gives that $-\ul X(e_q)$ has the exponential law with intensity $\Phi(q)$ and
\begin{equation}\label{eq:sup-lst}
        \ell_q(s):=< e^{-s\ol X(e_q)}>=\frac{q(s-\Phi(q))}{\Phi(q)(\psi(s)-q)}.
\end{equation}
The expression has the usual removable singularity at $s=\Phi(q)$.

\subsection{Embedded chain}

The unreflected process with resetting evolves as follows.  Immediately after a resetting Poissonian epoch $T_n$ with intensity $q$, the level is reset to $F$ of the level just before that epoch.
During the following inter-review interval, it evolves by adding an independent L\'evy increment.  We denote this L\'evy process with resetting by $Z(t)$. Let
\[
        Z_n:=Z(T_n-),
\]
be the pre-resetting level.  Then
\begin{equation}\label{eq:unref-chain}
        Z_{n+1}=F_{n+1}(Z_n)+\xi_{n+1},
        \qquad \xi_{n+1}\eqd X(e_q),
\end{equation}
with $\xi_{n+1}$ independent of $F_{n+1}$ and $Z_n$.  It is so-called one step analysis
based on observing the process at resetting times.
Assume that stationary law exists and is unique. Then, at stationarity,
\begin{equation}\label{eq:unref-fixed-point}
        Z\eqd F(Z)+\xi,
        \qquad \xi\eqd X(e_q),
\end{equation}
where $Z$ is a steady-state and it is independent of $F$ and $\xi$ on the right-hand side.
This identity is a direct consequence of
so-called PASTA property.

\subsection{Ergodicity conditions}

We start from identifying the sufficient ergodicity conditions for the existence of the stationary distribution.
that follows directly from the iterated-random-functions theorem.
\begin{theorem}[Unreflected contraction criterion]\label{thm:unref-contraction}
Assume that there are i.i.d.\ nonnegative random variables $A_n$ such that
\[
        |F_n(x)-F_n(y)|\leq A_n |x-y|\qquad\hbox{for all }x,y,
\]
and that, for some $x_0$,
\begin{equation}\label{eq:unref-contract-cond}
        <\log A_1><0,
        \qquad
        <log^+|F_1(x_0)+\xi_1|><\infty.
\end{equation}
Then the chain \eqref{eq:unref-chain} has a unique stationary distribution.  Moreover, for any initial state $Z_0$, $Z_n$ converges in distribution to the unique solution of \eqref{eq:unref-fixed-point}.
\end{theorem}

\begin{proof}
Let $G_n(x)=F_n(x)+\xi_n$.  Then
\[
        |G_n(x)-G_n(y)|=|F_n(x)-F_n(y)|\leq A_n|x-y|.
\]
Thus $G_n$ is a random Lipschitz map with negative mean logarithmic Lipschitz constant.  The second condition in \eqref{eq:unref-contract-cond} is the usual integrability condition.  The standard backward iteration
\[
        G_1\circ G_2\circ\cdots\circ G_n(x_0)
\]
converges almost surely to a finite limit that is independent of $x_0$, and the forward chain converges weakly to the law of that limit.
\end{proof}

For affine reset maps, the condition is straightforward: the contraction depends only on the multiplicative coefficient, while the additive part only enters through the integrability condition.

\begin{corollary}[Affine and random-affine resets]\label{cor:unref-affine-stability}
Let
\[
        F_n(x)=C_n x+B_n,
\]
where $(C_n,B_n)$ are i.i.d.\ and independent of $\xi_n$.  If
\[
        <\log|C_1|><0,
        \qquad
        <\log^+|B_1+\xi_1|><\infty,
\]
then the unreflected reset chain has a unique stationary distribution.  In particular, if $F(x)=cx+b$ with $|c|<1$ and $\E\log^+|b+X(e_q)|<\infty$, then the stationary distribution exists and is unique.
\end{corollary}

\subsection{Laplace transform of steady-state and moments recursion}
Let
\[
        v(s):=< e^{-sZ}>,
        \qquad
        h_F(s):=< e^{-sF(Z)}>.
\]
By \eqref{eq:unref-fixed-point} and independence,
\begin{equation}\label{eq:unref-general-lst}
        v(s)=m(s)h_F(s)=\frac{q}{q-\psi(s)}\,< e^{-sF(Z)}>,
        \qquad s\in\mathcal D_q.
\end{equation}
Assume that the necessary moments exist.  Define
\[
        \nu_n:=<\xi^n>=<X(e_q)^n>,
        \qquad
        m_n:=< Z^n>,
        \qquad
        \eta_n:=< F(Z)^n>.
\]
Expanding the fixed-point equation $Z\eqd \xi+F(Z)$ gives
\begin{equation}\label{eq:unref-moment-general}
        m_n=\sum_{k=0}^n \binom{n}{k}\nu_{n-k}\eta_k,
        \qquad n\geq0,
\end{equation}
with $m_0=\eta_0=\nu_0=1$.

If $F$ is polynomial or analytic and powers of $F(Z)$ can be expressed through the moments of $Z$, then \eqref{eq:unref-moment-general} is a recursive coefficient scheme.

\subsection{Affine reset: product transform and perpetuity representation}\label{subsec:unref-affine}

Let
\[
        F(x)=cx+b,
        \qquad |c|<1,
\]
be an affine resetting mechanism based on multiplying by fixed proprtion $c\in (-1,1)$ and shifting by $b\in \mathbb{R}$.
Then \eqref{eq:unref-general-lst} becomes
\begin{equation}\label{eq:unref-affine-functional}
        v(s)=\frac{q}{q-\psi(s)}e^{-bs}v(cs).
\end{equation}
Iterating \eqref{eq:unref-affine-functional} gives
\begin{equation}\label{eq:unref-affine-product}
        v(s)=\exp\left\{-\frac{bs}{1-c}\right\}
        \prod_{j=0}^{\infty}\frac{q}{q-\psi(c^j s)}.
\end{equation}
Equivalently,
\begin{equation}\label{eq:unref-affine-perpetuity}
        Z\eqd \frac{b}{1-c}+\sum_{j=0}^{\infty}c^j\xi_j,
\end{equation}
where $\xi_0,\xi_1,\ldots$ are i.i.d.\ copies of $X(e_q)$.  
In this case, when $F(x)=cx+b$, we also have
\[
        \eta_k=\E(cZ+b)^k=\sum_{r=0}^k \binom{k}{r}c^r b^{k-r}m_r.
\]
Thus, using \eqref{eq:unref-moment-general}, the term containing $m_n$ appears only when $k=n$ and $r=n$ and
\begin{equation}\label{eq:unref-affine-moment-recursion}
        (1-c^n)m_n
        =\sum_{k=0}^n \binom{n}{k}\nu_{n-k}
          \sum_{r=0}^{\min(k,n-1)}\binom{k}{r}c^r b^{k-r}m_r,
        \qquad n\geq1.
\end{equation}
This gives moments in increasing order.  The first two are
\begin{align}\label{eq:unref-affine-first-two}
        m_1&=\frac{b+\nu_1}{1-c},\\
        \Var (Z)&=\frac{\Var (\xi)}{1-c^2}.
\end{align}
%

\subsection{Unreflected Brownian under affine resetting}\label{subsec:unref-bm}
In this section we focus on very special case of the Brownian motion with affine resetting, that is,
Let
\[
        X(t)=dt+\sigma B(t),
        \qquad d\in\R,
        \quad \sigma>0,
\]
where $B$ is standard Brownian motion.  Then
\[
        \psi(s)=-ds+\frac{\sigma^2s^2}{2},
        \qquad
        m(s)=\frac{q}{q+ds-\sigma^2s^2/2}.
\]
For $F(x)=cx+b$, $|c|<1$, \eqref{eq:unref-affine-product} becomes
\begin{equation}\label{eq:unref-bm-product}
        v(s)=\exp\left\{-\frac{bs}{1-c}\right\}
        \prod_{j=0}^{\infty}
        \frac{q}{q+d c^j s-\frac{1}{2}\sigma^2c^{2j}s^2}.
\end{equation}
The product is valid for $s$ in a neighbourhood of the origin contained in the intersection of the domains of all factors.
The increment $\xi=X(e_q)$ has
\[
        <\xi>=\frac{d}{q},
        \qquad
        \Var(\xi)=\frac{\sigma^2}{q}+\frac{d^2}{q^2}.
\]
Consequently,
\begin{equation}\label{eq:unref-bm-mean-var}
        < Z>=\frac{b+d/q}{1-c},
        \qquad
        \Var (Z)=\frac{\sigma^2/q+d^2/q^2}{1-c^2}.
\end{equation}

\subsection{Unreflected finite-activity L\'evy flight under affine resetting}\label{subsec:unref-levyflight}

Here ``L\'evy flight'' is used in the finite-activity sense: jumps occur at Poisson times and have a regularly varying displacement distribution.  Let
\begin{equation}\label{eq:levy-flight-def}
        X(t)=at+\sum_{i=1}^{N_\lambda(t)}J_i,
\end{equation}
where $N_\lambda$ is Poisson with rate $\lambda$, $J_1,J_2,\ldots$ are i.i.d., and $a\in\R$.  Let
\[
        \hat G(s):=< e^{-sJ_1}>.
\]
Then
\begin{equation}\label{eq:levy-flight-laplace-exponent}
        \psi(s)=-as+\lambda(\hat G(s)-1),
        \qquad
        m(s)=\frac{q}{q+as-\lambda(\hat G(s)-1)}.
\end{equation}
For $F(x)=cx+b$, $|c|<1$,
\begin{equation}\label{eq:levy-flight-affine-product}
        v(s)=\exp\left\{-\frac{bs}{1-c}\right\}
        \prod_{j=0}^{\infty}
        \frac{q}{q+a c^j s-\lambda(\hat G(c^j s)-1)}.
\end{equation}

Assume $J$ has finite first two moments.  Then
\[
        <\xi>=\frac{a+\lambda< J>}{q},
\]
and, by conditioning on $e_q$,
\begin{align}\label{eq:levy-flight-increment-variance}
        \Var(\xi)
        &=\E\bigl[\Var(X(e_q)\mid e_q)\bigr]+\Var\bigl(<X(e_q)\mid e_q>\bigr)\notag\\
        &=\frac{\lambda< J^2>}{q}+\frac{(a+\lambda< J>)^2}{q^2}.
\end{align}
Therefore
\begin{equation}\label{eq:levy-flight-mean-var}
        < Z>=\frac{b+(a+\lambda< J>)/q}{1-c},
        \qquad
        \Var (Z)=\frac{\lambda< J^2>/q+(a+\lambda< J>)^2/q^2}{1-c^2}.
\end{equation}

If $J$ is nonnegative and regularly varying with index $-\alpha$, $\alpha>0$, i.e.
\[
        \Pp(J>x)=x^{-\alpha}L(x),
\]
where $L$ is slowly varying, then $\xi=X(e_q)$ has the one-big-jump tail
\begin{equation}\label{eq:levy-flight-increment-tail}
        \Pp(\xi>x)\sim \frac{\lambda}{q}\Pp(J>x).
\end{equation}
Using the perpetuity representation \eqref{eq:unref-affine-perpetuity} and standard subexponential closure,
\begin{equation}\label{eq:levy-flight-affine-tail}
        \Pp(Z>x)
        \sim \sum_{j=0}^{\infty}\Pp(c^j\xi_j>x)
        \sim \frac{\lambda}{q}\frac{1}{1-c^\alpha}\Pp(J>x),
        \qquad x\to\infty,
\end{equation}
for $0<c<1$.

\subsection{Unreflected spectrally-positive $\alpha$-stable flight under affine resetting}\label{subsec:unref-alpha-stable}

Let $X$ be a spectrally-positive stable L\'evy flight with index $\alpha\in(1,2)$ and negative linear component, normalised so that
\begin{equation}\label{eq:alpha-stable-exponent}
        \psi(s)=\mu s+\kappa s^\alpha,
        \qquad \mu>0,
        \quad \kappa>0.
\end{equation}
The positive jumps have tail
\begin{equation}\label{eq:alpha-stable-levy-tail}
        \Pi(x,\infty)\sim K_\alpha x^{-\alpha},
\end{equation}
for a constant $K_\alpha>0$ determined by the normalisation of the stable flight.

For $F(x)=cx+b$, $0<c<1$, the transform is
\begin{equation}\label{eq:alpha-stable-product}
        v(s)=\exp\left\{-\frac{bs}{1-c}\right\}
        \prod_{j=0}^{\infty}
        \frac{q}{q-\mu c^j s-\kappa c^{\alpha j}s^\alpha}.
\end{equation}
The domain is $0\leq s<\Phi(q)$, where $\Phi(q)$ is the positive root of $\mu s+\kappa s^\alpha=q$.
Since the positive jumps are regularly varying, the exponential-time increment satisfies
\begin{equation}\label{eq:alpha-stable-increment-tail}
        \Pp(X(e_q)>x)\sim \frac{1}{q}\Pi(x,\infty)
        \sim \frac{K_\alpha}{q}x^{-\alpha}.
\end{equation}
Therefore
\begin{equation}\label{eq:alpha-stable-Z-tail}
        \Pp(Z>x)\sim \frac{K_\alpha}{q(1-c^\alpha)}x^{-\alpha}.
\end{equation}


\subsection{Unreflected two-point random coefficient}\label{subsec:unref-two-point}

Let
\begin{equation}\label{eq:two-point-C}
        F(x)=Cx+b,
        \qquad
        \Pp(C=c)=p,
        \qquad
        \Pp(C=0)=1-p,
\end{equation}
with $0\leq c<1$, $0\leq p\leq1$, and $b\geq0$.  Then
\begin{equation}\label{eq:two-point-functional}
        v(s)=m(s)e^{-bs}\bigl((1-p)+p v(cs)\bigr).
\end{equation}
Iterating gives the explicit series
\begin{equation}\label{eq:two-point-series}
        v(s)=(1-p)\sum_{n=0}^{\infty}p^n
        \exp\left\{-b s\frac{1-c^{n+1}}{1-c}\right\}
        \prod_{j=0}^{n}m(c^j s).
\end{equation}
For $p=1$, the geometric termination disappears and \eqref{eq:two-point-series} reduces to the affine product \eqref{eq:unref-affine-product}.

The corresponding perpetuity representation is
\begin{equation}\label{eq:two-point-perpetuity}
        Z\eqd \sum_{j=0}^{N}\left(\prod_{i=0}^{j-1}C_i\right)(b+\xi_j),
\end{equation}
where $N$ is geometric with
\[
        \Pp(N=n)=(1-p)p^n,
        \qquad n=0,1,2,\ldots,
\]
and, on $\{N=n\}$, $C_0=\cdots=C_{n-1}=c$ and $C_n=0$.

For moments, let $\gamma_r:=< C^r>=p c^r$ for $r\geq1$, with $\gamma_0=1$.  Then
\begin{equation}\label{eq:random-affine-moment-recursion}
        (1-\gamma_n)m_n
        =\sum_{k=0}^{n}\binom{n}{k}\nu_{n-k}
        \sum_{r=0}^{\min(k,n-1)}\binom{k}{r}\E(C^r b^{k-r})m_r,
        \qquad n\geq1.
\end{equation}
Because $b$ is deterministic here,
\[
        <C^r b^{k-r}>=b^{k-r}\gamma_r.
\]
The same recursion works for any finite-support $C$.  For a three-point reset coefficient
\[
        \Pp(C=c_i)=p_i,
        \qquad i=1,2,3,
\]
one only replaces $\gamma_r$ by $p_1c_1^r+p_2c_2^r+p_3c_3^r$.  The transform equation becomes
\[
        v(s)=m(s)e^{-bs}\sum_{i=1}^3p_i v(c_i s),
\]
which is a branching functional equation.

If $\xi$ is regularly varying with index $-\alpha$ and $B\equiv b$ is light-tailed, then
\begin{equation}\label{eq:random-affine-tail}
        \Pp(Z>x)\sim \frac{1}{1-< C^\alpha>}\Pp(\xi>x).
\end{equation}
For the two-point coefficient \eqref{eq:two-point-C}, this gives
\begin{equation}\label{eq:two-point-tail}
        \Pp(Z>x)\sim \frac{1}{1-pc^\alpha}\Pp(\xi>x).
\end{equation}

\section{The reflected spectrally positive reset model}\label{sec:reflected}

\subsection{Model and embedded chain}

Assume now that $X$ is spectrally positive and not a subordinator.  Given an initial level $x\geq0$, define the reflected process at $0$ by
\begin{equation}\label{eq:reflected-process}
        Y(t)=x+X(t)-\inf_{0\leq u\leq t}(x+X(u))\wedge0.
\end{equation}
At resetting Poissonian epochs with intensity $q>0$, the resetting using functional $F$ is performed.  Let
\[
        Y_n:=Y(T_n-).
\]
Then
\begin{equation}\label{eq:reflected-embedded}
        Y_{n+1}\eqd \max\{F_{n+1}(Y_n)+X(e_q),\ol X(e_q)\}.
\end{equation}
and at stationarity (assuming that it exists and is unique),
\begin{equation}\label{eq:reflected-fixed-point}
        Y\eqd \max\{F(Y)+X(e_q),\ol X(e_q)\},
\end{equation}
where by $Y$ we denote the steady-state.
Equivalently, by Wiener--Hopf factorisation,
\begin{equation}\label{eq:BKP-style-recursion}
        Y\eqd \ol X(e_q)+\bigl(F(Y)-E_\Phi\bigr)^+,
\end{equation}
where $E_\Phi$ has the exponential distribution with parameter $\Phi(q)$ and $\ol X(e_q)$, $E_\Phi$, and $F(Y)$ are independent on the right-hand side.

\subsection{Ergodicity conditions}

The following is the Palmowski-Vlasiou \cite{PalmowskiVlasiou2011} comparison criterion, written in the notation used here.

\begin{theorem}[Reflected comparison criterion]\label{thm:reflected-stability}
Assume that $X$ is not a subordinator and that there exist nonnegative i.i.d.\ random variables $(A_n,B_n)$ such that
\begin{equation}\label{eq:reflected-domination}
        F_n(x)\leq A_nx+B_n,
        \qquad x\geq0.
\end{equation}
Then the reflected reset chain has a unique stationary limiting distribution if one of the following holds:
\begin{enumerate}[label=(\alph*)]
\item $A_n=0$ a.s.;
\item $A_n=1$ a.s. and $< B_1>+< X(1)>/q<0$;
\item $<\log A_1><0$ and $<\log^+B_1><\infty$.
\end{enumerate}
Here, if $B_1=0$ in (c), one may dominate by $A_nx+1$ without changing the conclusion.
\end{theorem}

\begin{proof}[Proof sketch]
Under \eqref{eq:reflected-domination}, the original chain is stochastically dominated by the affine reflected chain
\[
        Y_{n+1}^{a}\eqd \max\{A_{n+1}Y_n^a+B_{n+1}+X(e_q),\ol X(e_q)\}.
\]
Cases (a)--(c) imply return to a compact set.  Since $X$ is not a subordinator, the reflected process can hit zero, and zero is a regeneration state in the Harris-chain sense.  The minorisation on a compact set then gives positive Harris recurrence and uniqueness of the stationary distribution.
\end{proof}

\begin{remark}\rm
For the partial restting $F(x)=Cx$ with $0\leq C\leq1$ and $\Pp(C<1)>0$, condition (c) holds after domination by $Cx+1$ whenever $<\log C><0$.
Note that this holds for example when $C$ has a uniform law $U(0,1)$, because in this case $<\log C>=-1$.
\end{remark}

\subsection{Laplace transform, moments and tail of steady-state}
Let
\[
        v(s):=< e^{-sY}>,
        \qquad
        h_F(s):=< e^{-sF(Y)}>.
\]
Using \eqref{eq:BKP-style-recursion}, independence, and the elementary identity
\begin{equation}\label{eq:exp-trunc-transform}
        < e^{-s(a-E_\Phi)^+}>
        =\frac{\Phi(q)e^{-sa}-s e^{-\Phi(q)a}}{\Phi(q)-s},
        \qquad a\geq0,
\end{equation}
we get
\begin{align}\label{eq:reflected-general-lst-long}
        v(s)
        &=\ell_q(s)\,\frac{\Phi(q)h_F(s)-s h_F(\Phi(q))}{\Phi(q)-s},
\end{align}
where $\ell_q(s)$ is defined in \eqref{eq:sup-lst}.
Substituting \eqref{eq:sup-lst} and simplifying yields the reset-isolated form
\begin{equation}\label{eq:reflected-general-lst}
        v(s)=\frac{q}{\Phi(q)(q-\psi(s))}
        \left\{\Phi(q)h_F(s)-s h_F(\Phi(q))\right\}.
\end{equation}
Equivalently,
\begin{equation}\label{eq:reflected-general-lst-a}
        v(s)=a(s)\left(h_F(s)-\frac{s}{\Phi(q)}h_F(\Phi(q))\right),
        \qquad
        a(s):=\frac{q}{q-\psi(s)}.
\end{equation}

Let $\xi:=X(e_q)$ and write
\[
        \nu_n:=<\xi^n,>
        \qquad
        m_n:=< Y^n>,
        \qquad
        \eta_n:=< F(Y)^n>,
        \qquad
        K_F:=h_F(\Phi(q))=< e^{-\Phi(q)F(Y)}>.
\]
Expanding \eqref{eq:reflected-general-lst-a} gives, for $n\geq1$,
\begin{equation}\label{eq:reflected-moment-general}
        m_n=\sum_{k=0}^n \binom{n}{k}\nu_{n-k}\eta_k
        +\frac{nK_F}{\Phi(q)}\nu_{n-1}.
\end{equation}
This is not yet closed for a general $F$, but it becomes recursive for affine or random-affine reset maps.  For $F(x)=cx$, $0\leq c<1$, $\eta_k=c^k m_k$ and
\begin{equation}\label{eq:reflected-cx-moment-recursion}
        (1-c^n)m_n
        =\sum_{k=0}^{n-1}\binom{n}{k}\nu_{n-k}c^k m_k
        +\frac{nK_c}{\Phi(q)}\nu_{n-1},
        \qquad
        K_c:=v(c\Phi(q)).
\end{equation}
In particular,
\begin{equation}\label{eq:reflected-cx-mean}
        < Y>=\frac{\nu_1+K_c/\Phi(q)}{1-c}
        =\frac{< X(e_q)>+v(c\Phi(q))/\Phi(q)}{1-c}.
\end{equation}
Since $\nu_1=< X(1)>/q=-\psi'(0)/q$, this may also be written as
\[
        < Y>=\frac{v(c\Phi(q))/\Phi(q)-\psi'(0)/q}{1-c}.
\]
For $F(x)=Cx$ with $C$ independent of $Y$,
\begin{equation}\label{eq:reflected-randomC-moment-recursion}
        (1-< C^n>)m_n
        =\sum_{k=0}^{n-1}\binom{n}{k}\nu_{n-k}< C^k>\,m_k
        +\frac{nK_C}{\Phi(q)}\nu_{n-1},
        \qquad
        K_C:=< v(C\Phi(q))>.
\end{equation}
For $C$ having a uniform distribution $U(0,1)$, we have $< C^n>=(n+1)^{-1}$ and $K_C=\Phi(q)^{-1}\int_0^{\Phi(q)}v(u)\,\dd u$.

Following Palmowski--Vlasiou \cite{PalmowskiVlasiou2011}, 
we record the regularly varying tail expression for the most relevant partial resetting case with random proportion.
\begin{proposition}[Regularly varying reflected L\'evy process under partial resetting]\label{prop:reflected-proportional-tail}
Assume that $X$ is spectrally positive, not a subordinator, and that its L\'evy measure satisfies
\[
        \ol\Pi(x):=\Pi(x,\infty)\in \RV_{-\alpha},
        \qquad \alpha>0.
\]
If
$F(x)=Cx$, where $0\leq C\leq1$,
then
\begin{equation}\label{eq:reflected-proportional-tail}
        \Pp(Y>x)
        \sim \frac{1}{q(1-
< C^\alpha>)}\ol\Pi(x).
\end{equation}
\end{proposition}

\begin{proof}[Proof sketch]
The Wiener--Hopf representation \eqref{eq:BKP-style-recursion} gives
\[
        Y\eqd \ol X(e_q)+(CY-E_\Phi)^+.
\]
The exponential $E_\Phi$ is light-tailed and does not change a regularly varying tail.  The stochastic recurrence is therefore tail-equivalent to the random-coefficient affine recursion
\[
        Y\approx \ol X(e_q)+CY.
\]
For regularly varying additive term with index $\alpha$ and $< C^\alpha><1$, the standard one-large-term argument gives
\[
        \Pp(Z>x)\sim \sum_{j=0}^{\infty}\Pp(C_1\cdots C_j\,\ol X_j(e_q)>x)
        \sim \sum_{j=0}^{\infty}(< C^\alpha>)^j\Pp(\ol X(e_q)>x).
\]
We also know, by principle of one big jump, that
\begin{equation}\label{eq:sup-tail-onebigjump}
        \Pp(\ol X(e_q)>x)\sim \frac{1}{q}\ol\Pi(x).
\end{equation}
Using \eqref{eq:sup-tail-onebigjump} gives \eqref{eq:reflected-proportional-tail}.
\end{proof}
For deterministic $C=c\in[0,1)$ above proposition gives
\begin{equation}\label{eq:reflected-deterministic-c-tail}
        \Pp(Y>x)
        \sim \frac{1}{q(1-c^\alpha)}\ol\Pi(x),
\end{equation}
and for $C$ having uniform law $U(0,1)$,
\begin{equation}\label{eq:reflected-uniform-tail}
        \Pp(Y>x)
        \sim \frac{\alpha+1}{\alpha q}\ol\Pi(x).
\end{equation}


\subsection{The partial resetting 
}\label{subsec:reflected-cx}
Let $F(x)=cx$, $0\leq c<1$.  Then $h_F(s)=v(cs)$ and \eqref{eq:reflected-general-lst-a} gives
\begin{equation}\label{eq:reflected-cx-functional}
        v(s)=a(s)\left(v(cs)-\frac{s}{\Phi(q)}v(c\Phi(q))\right),
        \qquad
        a(s)=\frac{q}{q-\psi(s)}.
\end{equation}
Let
\[
        K:=v(c\Phi(q)),
        \qquad
        A_\infty(s):=\prod_{j=0}^{\infty}a(c^js),
\]
and
\[
        B_\infty(s):=\sum_{k=0}^{\infty}c^k\prod_{j=0}^{k}a(c^js).
\]
Iterating \eqref{eq:reflected-cx-functional} yields
\begin{equation}\label{eq:reflected-cx-product-solution}
        v(s)=A_\infty(s)-K\frac{s}{\Phi(q)}B_\infty(s).
\end{equation}
The constant $K$ is obtained by setting $s=c\Phi(q)$:
\begin{equation}\label{eq:reflected-cx-K}
        K=\frac{A_\infty(c\Phi(q))}{1+cB_\infty(c\Phi(q))}.
\end{equation}
Equations \eqref{eq:reflected-cx-product-solution}--\eqref{eq:reflected-cx-K} give a completely explicit transform in terms of the L\'evy exponent $\psi$.

For completeness, again using Palmowski--Vlasiou \cite{PalmowskiVlasiou2011}, we record the exponential-type asymptotics of the tail of the law of $Y$.  Suppose there exists $\theta_*>0$ such that
\begin{equation}\label{eq:cramer-root}
        \psi(\theta_*)=q,
\end{equation}
and the corresponding tilted ladder-height mean $m_q$ is finite.
Let \begin{equation}\label{asslighttail}
        < e^{\theta_*F(Y)}><\infty.
\end{equation}
Then the reflected L\'evty process under resetting the asymptotic tail is given by
\begin{equation}\label{eq:cramer-tail}
        \Pp(Y>x)\sim D e^{-\theta_*x},
\end{equation}
with
\begin{equation}\label{eq:cramer-constant}
        D=\Pp(E_\Phi>F(Z))\frac{\kappa(q,0)}{\theta_*m_q},
        \qquad
        E_\Phi\sim\Exp(\Phi(q)),
\end{equation}
where $\kappa(q,0)$ is LST of ladder time.
For bounded resets and for $F(x)=(B-x)^+$ with light-tailed $B$, the condition \eqref{asslighttail}
is immediate.  For $F(x)=Cx$, it becomes a condition on $v(-\theta_*C)$, and therefore requires finite positive exponential moments of $Y$.

\subsection{Reflected Brownian motion under partial resettting}\label{subsec:reflected-bm-cx}
Consider $F(x)=cx$.
Let
\[
        X(t)=\sigma B(t)-\mu t,
        \qquad \mu>0,
        \quad \sigma>0.
\]
Then
\[
        \psi(s)=\mu s+\frac{\sigma^2s^2}{2},
        \qquad
        \Phi(q)=\frac{-\mu+\sqrt{\mu^2+2q\sigma^2}}{\sigma^2}.
\]
Thus
\begin{equation}\label{eq:reflected-bm-cx-a}
        a(s)=\frac{q}{q-\mu s-\sigma^2s^2/2}.
\end{equation}
Substituting this $a$ in \eqref{eq:reflected-cx-product-solution} gives
\begin{equation}\label{eq:reflected-bm-cx-transform}
        v(s)=\prod_{j=0}^{\infty}\frac{q}{q-\mu c^j s-\frac{1}{2}\sigma^2c^{2j}s^2}
        -K\frac{s}{\Phi(q)}
        \sum_{k=0}^{\infty}c^k
        \prod_{j=0}^{k}\frac{q}{q-\mu c^j s-\frac{1}{2}\sigma^2c^{2j}s^2},
\end{equation}
where $K$ is given by \eqref{eq:reflected-cx-K}.  The first moment follows from \eqref{eq:reflected-cx-mean}:
\begin{equation}\label{eq:reflected-bm-cx-mean}
        < Y>=\frac{K/\Phi(q)-\mu/q}{1-c}.
\end{equation}
The second moment follows from \eqref{eq:reflected-cx-moment-recursion}.  Since
\[
        \nu_1=< X(e_q)>=-\frac{\mu}{q},
        \qquad
        \nu_2=< X(e_q)^2>=\frac{\sigma^2}{q}+\frac{2\mu^2}{q^2},
\]
we have
\begin{equation}\label{eq:reflected-bm-cx-second-moment}
        (1-c^2)m_2
        =\nu_2+2c\nu_1m_1+\frac{2K}{\Phi(q)}\nu_1.
\end{equation}
Together with \eqref{eq:reflected-bm-cx-mean}, this gives $\Var (Y)=m_2-m_1^2$.

\subsection{Reflected reset $F(x)=(B-x)^+$ with exponential $B$}\label{subsec:reflected-Bminusx}

Let $F(x)=(B-x)^+$, where $B$ has an exponential law $\Exp(\beta)$ that is independent of $Y$.  The memoryless property closes the reset transform.  Conditional on $Y=z$,
\[
        (B-z)^+=0 \quad\hbox{on }\{B\leq z\},
\]
and, on $\{B>z\}$, the excess $B-z$ is exponential with rate $\beta$.  Therefore
\begin{equation}\label{eq:BminusZ-tail}
        \Pp((B-Y)^+>y)=< e^{-\beta(Y+y)}>=v(\beta)e^{-\beta y},
        \qquad y\geq0.
\end{equation}
Thus $(B-Y)^+$ has an atom $1-v(\beta)$ at zero and an exponential density $v(\beta)\beta e^{-\beta y}$ on $(0,\infty)$.  Hence
\begin{equation}\label{eq:BminusZ-transform}
        h_F(s)=< e^{-s(B-Z)^+}>
        =1-v(\beta)+v(\beta)\frac{\beta}{\beta+s}
        =1-v(\beta)\frac{s}{\beta+s}.
\end{equation}
Let $r:=v(\beta)$ and $\Phi:=\Phi(q)$.  Substituting \eqref{eq:BminusZ-transform} into \eqref{eq:reflected-general-lst} gives
\begin{equation}\label{eq:BminusZ-vs}
        v(s)=\frac{q}{\Phi(q-\psi(s))}
        \left[\Phi\left(1-r\frac{s}{\beta+s}\right)-s\left(1-r\frac{\Phi}{\beta+\Phi}\right)\right].
\end{equation}
The remaining scalar $r=v(\beta)$ is obtained by setting $s=\beta$ in \eqref{eq:BminusZ-vs}.  With
\[
        D_\beta:=\frac{q}{\Phi(q-\psi(\beta))},
\]
we get
\begin{align*}
        r
        &=D_\beta\left[\Phi\left(1-\frac{r}{2}\right)
          -\beta\left(1-r\frac{\Phi}{\beta+\Phi}\right)\right]\\
        &=D_\beta(\Phi-\beta)-D_\beta r\,\Phi\left(\frac{1}{2}-\frac{\beta}{\beta+\Phi}\right).
\end{align*}
Thus, provided $\beta<\Phi$ and $\beta$ is in the transform domain,
\begin{equation}\label{eq:BminusZ-r}
        r=\frac{D_\beta(\Phi-\beta)}{1+D_\beta\Phi(\Phi-\beta)/(2(\beta+\Phi))}.
\end{equation}
Equations \eqref{eq:BminusZ-vs} and \eqref{eq:BminusZ-r} provide the LST explicitly.

\subsection{Partial resetting with uniform proportion}\label{subsec:reflected-uniform}
Let $F(x)=Cx$ with $C$ having uniform distribution $U(0,1)$.
Then for $C$, independent of $Y$ we have
\begin{equation}\label{eq:uniform-hF}
        h_F(s)=< e^{-sCY}>=\int_0^1 v(su)\,\dd u=\frac{1}{s}\int_0^s v(y)\,\dd y.
\end{equation}
Write
\[
        F_0(s):=\int_0^s v(y)\,\dd y,
        \qquad
        \Phi:=\Phi(q).
\]
Then \eqref{eq:reflected-general-lst-a} becomes
\begin{equation}\label{eq:uniform-diffeq-raw}
        v(s)=\frac{F_0(s)/s-(s/\Phi^2)F_0(\Phi)}{1-\psi(s)/q}.
\end{equation}
Since $v=F_0'$, this is the first-order inhomogeneous differential equation
\begin{equation}\label{eq:uniform-diffeq}
        F_0'(s)=\frac{F_0(s)-b_0s^2/2}{s(1-\psi(s)/q)},
        \qquad
        b_0:=\frac{2F_0(\Phi)}{\Phi^2}.
\end{equation}
Solving \eqref{eq:uniform-diffeq} with the boundary condition $F_0(0)=0$ and the regularity condition at $s=\Phi$ gives 
final form.  Define
\begin{equation}\label{eq:uniform-g}
        g(s):=\int_0^s
        \exp\left\{-\int_0^x \frac{\psi(y)}{y(q-\psi(y))}\,\dd y\right\}\dd x.
\end{equation}
Then
\begin{equation}\label{eq:uniform-solution}
        v(s)=\frac{1-g(s)/g(\Phi)}{g'(s)(1-\psi(s)/q)},
        \qquad 0\leq s<\Phi.
\end{equation}
Moreover $b_0=1/g(\Phi)$.

\subsection{Reflected Brownian motion under uniform partial resetting}
Take again
\[
        X(t)=\sigma B(t)-\mu t,
        \qquad
        \psi(s)=\mu s+\frac{\sigma^2s^2}{2}.
\]
Let
\begin{equation}\label{eq:uniform-bm-roots}
        y_1:=\Phi(q)=\frac{-\mu+\sqrt{\mu^2+2q\sigma^2}}{\sigma^2},
        \qquad
        y_2:=\frac{-\mu-\sqrt{\mu^2+2q\sigma^2}}{\sigma^2}<0.
\end{equation}
Then
\begin{equation}\label{eq:uniform-bm-partial-fraction}
        -\frac{\psi(y)}{y(q-\psi(y))}
        =\frac{D_1}{y-y_1}+\frac{D_2}{y-y_2},
        \qquad
        D_1=\frac{-y_2}{y_1-y_2},
        \quad
        D_2=\frac{y_1}{y_1-y_2}.
\end{equation}
Hence
\begin{equation}\label{eq:uniform-bm-g}
        g(s)=\int_0^s
        \left(\frac{x-y_1}{-y_1}\right)^{D_1}
        \left(\frac{x-y_2}{-y_2}\right)^{D_2}\dd x,
        \qquad 0\leq s\leq y_1.
\end{equation}
Substitution into \eqref{eq:uniform-solution} gives
\begin{equation}\label{eq:uniform-bm-transform}
        v(s)=
        \left(\frac{y_1}{y_1-s}\right)^{1+D_1}
        \left(\frac{-y_2}{-y_2+s}\right)^{1+D_2}
        \left(1-\frac{g(s)}{g(y_1)}\right),
        \qquad 0\leq s<y_1.
\end{equation}
The integral $g$ may be written through incomplete beta functions.  With
\[
        z(s)=\frac{s-y_2}{y_1-y_2},
\]
we have
\begin{equation}\label{eq:uniform-bm-incomplete-beta}
\begin{aligned}
        \int_0^s (x-y_1)^{D_1}(x-y_2)^{D_2}\,\dd x
        &=(-1)^{D_1}(y_1-y_2)^{D_1+D_2+1}\\
        &\quad\times\left[ B\bigl(z(s);1+D_2,1+D_1\bigr)
          -B\bigl(z(0);1+D_2,1+D_1\bigr)\right],
\end{aligned}
\end{equation}
where $B(z;a,b)=\int_0^z t^{a-1}(1-t)^{b-1}\,\dd t$.

%

\subsection{Reflected spectrally positive stable flight under partial restting}\label{subsec:reflected-stable-collapse}

Let $X$ have exponent
\[
        \psi(s)=\mu s+\kappa s^\alpha,
        \qquad \alpha\in(1,2),
\]
as in \eqref{eq:alpha-stable-exponent}.  For deterministic proportion of the partial resetting, hence for $F(x)=cx$, $0<c<1$, the transform is given by \eqref{eq:reflected-cx-product-solution} with
\begin{equation}\label{eq:reflected-stable-a}
        a(s)=\frac{q}{q-\mu s-\kappa s^\alpha}
\end{equation}
and $\Phi(q)$ solving $\mu\Phi(q)+\kappa\Phi(q)^\alpha=q$. 
Note that the L\'evy tail satisfies \eqref{eq:alpha-stable-levy-tail}. Then Proposition \ref{prop:reflected-proportional-tail} gives
\begin{equation}\label{eq:reflected-stable-tail}
        \Pp(Y>x)
        \sim \frac{K_\alpha}{q(1-c^\alpha)}x^{-\alpha}.
\end{equation}
For uniform partial resetting, one should replace $(1-c^\alpha)^{-1}$ by $(\alpha+1)/\alpha$.

\section{Conclusion}
The unreflected model with resetting gives a clean baseline: one step analysis reduce the description of stationary state $Z$ to the fixed-point equation $Z\eqd X(e_q)+F(Z)$ for the steady-state $Z$.  The reflected model keeps the same reset structure but inserts the Wiener--Hopf correction term, producing the stationary state identity \eqref{eq:BKP-style-recursion}.  
The presented method is quite robust. For example, we can analyse other processes under resetting:
unreflected L\'evy process under affine random resetting, that is, $F(x)=Cx+B$ with $C$ having finite support, beta distributions, or mixtures with an atom at zero.
Similarly, we can handle polynomial resetting mechanism with $F(x)=b_0+b_1x+\cdots+b_dx^d$. Then \eqref{eq:unref-moment-general} expresses $m_n$ through moments up to order $dn$.
Very interesting seems to be the case the resetting is based on finite-support random functions.  In this case $F$ takes maps $f_1,\ldots,f_r$ with probabilities $p_1,\ldots,p_r$.
Then
\[
        v(s)=m(s)\sum_{i=1}^r p_i< e^{-sf_i(Z)}>.
\]
The same resetting mechanisms can be considered for the reflected L\'evy process.
Within this context quite intriguing is the case of Brownian plus compound Poisson with phase-type positive jumps which yields rational $\psi$.
Then all transform formulae reduce to products or finite sums involving rational factors.  This is the natural setting for numerical inversion.

\end{document}